# A Bayesian test for excess zeros in a zero-inflated power series distribution[*]

## Archan Bhattacharya[1], Bertrand S. Clarke[2] and Gauri S. Datta[1]

*University of Georgia, University of British Columbia and University of Georgia*

**Abstract:** Power series distributions form a useful subclass of one-parameter discrete exponential families suitable for modeling count data. A zero-inflated power series distribution is a mixture of a power series distribution and a degenerate distribution at zero, with a mixing probability $p$ for the degenerate distribution. This distribution is useful for modeling count data that may have extra zeros. One question is whether the mixture model can be reduced to the power series portion, corresponding to $p = 0$, or whether there are so many zeros in the data that zero inflation relative to the pure power series distribution must be included in the model i.e., $p \geq 0$. The problem is difficult partially because $p = 0$ is a boundary point.

Here, we present a Bayesian test for this problem based on recognizing that the parameter space can be expanded to allow $p$ to be negative. Negative values of $p$ are inconsistent with the interpretation of $p$ as a mixing probability, however, they index distributions that are physically and probabilistically meaningful. We compare our Bayesian solution to two standard frequentist testing procedures and find that using a posterior probability as a test statistic has slightly higher power on the most important ranges of the sample size $n$ and parameter values than the score test and likelihood ratio test in simulations. Our method also performs well on three real data sets.

## 1. Zero-inflated families

Models for count data often fail to fit in practice because of the presence of more zeros in the data than is explained by a standard model. This situation is often called zero inflation because the number of zeros is inflated from the baseline number of zeros that would be expected in, say, a one-parameter discrete exponential family. Zero inflation is a special case of overdispersion that contradicts the relationship between the mean and variance in a one-parameter exponential family. One way to address this is to use a two-parameter distribution so that the extra parameter permits a larger variance. Efron [9] developed the notion of double exponential family, a two-parameter modification of a standard one-parameter exponential family, that allows a higher variance than permitted by the one-parameter version. This is reasonable in some examples, typical count data distributions, such as Poisson, cannot be used to model data containing extra zeros.

Johnson, Kotz and Kemp ([13], pages 312–318) discuss a simple modification of a power series (PS) distribution $f(\cdot|\theta)$ to handle extra zeros. An extra proportion of

---

[*]Supported in part by NSF Grant SES-02-41651.

[1]Department of Statistics, University of Georgia, Athens, GA 30602-1952, USA, e-mail: archanbh@gmail.com; gaurisdatta@gmail.com

[2]Department of Statistics, 333-6356 Agricultural Road, University of British Columbia, Vancouver, BC, V6T 1Z2 Canada, e-mail: riffraff@stat.ubc.ca

*AMS 2000 subject classifications:* Primary 62F15, 62F03; secondary 62F05.

*Keywords and phrases:* Jeffreys' prior, posterior distribution, score test.

zeros, $p$, is added to the proportion of zeros from the original discrete distribution, while decreasing the remaining proportions in an appropriate way. So the zero-inflated PS distribution is defined as

$$(1) \qquad f^*(y|p, \theta) = \begin{cases} p + (1-p)f(0|\theta), & \text{if } y = 0, \\ (1-p)f(y|\theta), & \text{if } y > 0, \end{cases}$$

where $\theta \in \Theta$, the parameter space and the mixing parameter $p$ ranges over the interval

$$-f(0|\theta)/(1 - f(0|\theta)) < p < 1.$$

This allows the distribution to be well defined for certain negative values of $p$, depending on $\theta$. Although the mixing interpretation is lost when $p < 0$, these values have a natural interpretation in terms of zero-deflation, relative to a PS model. Correspondingly, $p > 0$ can be regarded as zero inflation relative to a PS model. Note that the PS family contains all discrete one-parameter exponential families so an appropriate choice of PS model in (1) permits any desired interpretation for the data corresponding to the second term. The first term allows an extra proportion $p$ of zeros to be added to the discrete PS distribution; this data is effectively regarded as a sort of contamination. Note that, zero inflation (zero deflation, respectively) does not imply that model (1) has larger (smaller, respectively) variance than the non-inflated version.

The first question to be asked is whether the degenerate distribution at zero is necessary. If it is not, then no zero inflation needs to be modeled and the model simplifies to $f(y|\theta)$. Clearly, this is a hypothesis testing problem. If $p$ is not allowed to be negative, $p = 0$ is a boundary point and testing $\mathcal{H}_0 : p = 0$ vs. $\mathcal{H}_1 : p > 0$ is a notoriously difficult problem for both Bayesians and frequentists for which few results are available. (See Self and Liang [18] and Silvapulle and Silvapulle [19] for some asymptotics from a frequentist perspective.) Permitting negative values of $p$ removes the boundary point problem so that the analytic challenges become manageable. The Bayes test obtained here compares favorably with standard frequentist methods in the real and simulated data cases we have examined.

Familiar cases in which testing $\mathcal{H}_0 : p = 0$ is useful include the zero-inflated Poisson (ZIP) distribution with parameters $(p, \theta)$ given by

$$(2) \qquad f^*(y|p, \theta) = pI_{\{y=0\}} + (1-p)\frac{e^{-\theta}\theta^y}{y!}, \qquad y = 0, 1, 2, \ldots$$

in which $\theta > 0$, $\frac{-e^{-\theta}}{1 - e^{-\theta}} < p < 1$ and $E(Y|p, \theta) = (1-p)\theta$ and the zero-inflated geometric distribution with parameters $(p, \theta)$:

$$(3) \qquad f^*(y|p, \theta) = pI_{\{y=0\}} + (1-p)(1-\theta)\theta^y, \qquad y = 0, 1, 2, \ldots$$

in which $0 < \theta < 1$, $-\frac{1-\theta}{\theta} < p < 1$, and $E(Y|p, \theta) = (1-p)\theta/(1-\theta)$. The zero-inflated binomial is similar.

These models have been examined from a frequentist standpoint. The earliest results on zero inflation can be found in Cochran [4] and Rao and Chakravarti [17]. In fitting a Poisson model to count data these authors checked whether lack of fit was due to the presence of extra zeros in the data by using an exact test and likelihood ratio test. Also in the context of a ZIP model, El-Shaarawi [10] obtained the ML estimator and used its asymptotic distribution to construct a confidence interval for the mean parameter. A peculiarity of the MLE for $p$ is that it can give negative values if there are no zeros in the data. Van Broek [1] derived the score test for the zero inflation parameter $p$ for testing $\mathcal{H}_0 : p = 0$ vs. $\mathcal{H}_1 : p \neq 0$. The two-sided alternative, however, gives up some power because the desired alternative is one-sided $\mathcal{H}_1 : p > 0$. A secondary problem is that the performance of this test deteriorates as the mean parameter increases. This may not be a serious problem because, as the mean increases, excess zeros will become more visually obvious since the Poisson model assigns ever less probability to zero.

More generally, Deng and Paul [7] extended the score test to general one-parameter exponential family. Thus, motivated by industrial applications, they studied a regression model for the mean parameter of the exponential distribution. Later, Deng and Paul [8] treated overdispersion and zero inflation simultaneously. In the ZIP context, Lambert [14] fitted a logistic regression model for $p$ and a log-linear model for $\theta$, using an EM algorithm to obtain estimates. Hall [12] extended this approach by adding random effects to the ZIP model and considered the case of a zero-inflated binomial model as well.

From the Bayesian standpoint, Ghosh, Mukhopadhyay and Lu [11] estimated the parameters in a ZIP model in regression context as an alternative to traditionally used maximum likelihood based methods. Their simulation studies showed the Bayesian method had better finite sample performance than the classical method, giving tighter interval estimates and higher coverage probabilities. Our work can be regarded as a continuation of their work for hypothesis testing.

In this paper, the main goal is to give a Bayes test of $\mathcal{H}_0 : p = 0$ vs. $\mathcal{H}_1 : p > 0$. To this end we consider the posterior probability

$$(4) \qquad T(\mathbf{Y}) = P(p > 0 | \mathbf{Y}) = \frac{\int_\Theta \int_0^1 L(p, \theta) \pi(p|\theta) \pi(\theta) dp d\theta}{\int_\Theta \int_{\frac{-f(0|\theta)}{1-f(0|\theta)}}^1 L(p, \theta) \pi(p|\theta) \pi(\theta) dp d\theta},$$

in which $L(p, \theta)$ is the likelihood function from a zero-inflated PS model and $\mathbf{Y}$ is a vector of $n$ data points. The corresponding rejection region is $T(\mathbf{Y}) > c$ for some suitable $c$. Asymptotic choice of $c$ is discussed in Section 3. Using (4) necessitates careful consideration of prior selection so that neither the null nor the alternative hypothesis will be unduly favored. This is done here by using Jeffreys' prior.

Treating (4) as a frequentist test statistic, we derive some of its properties. In particular, we obtain higher order corrections for its asymptotic behavior. Then, we verify computationally that the power, a frequentist property, of the Bayes test for the ZIP family is roughly the same or a little higher than the power of the score test and the likelihood ratio test, for the hardest and most important ranges of $n$, $p$ and $\theta$ i.e., small to moderate $p$, small-ish $\theta$, and small to moderate $n$. This is striking because Jeffreys' priors favor small $\theta$'s and $p$'s near 0 and 1, and so are relatively unfavorable to the null. From the estimation standpoint, we verify that the posterior density is well behaved and gives reasonable credible intervals. As a final verification, we apply our techniques to three real data sets computing Bayes factors and score tests for the presence of zero inflation and obtaining estimates for the zero inflation as appropriate.

The structure of the paper is as follows. In Section 2 we provide some background on the properties of zero-inflated models from a Bayesian standpoint. In Section 3 we present the Bayesian test and give some of its properties. In Section 4 we develop Bayes estimation. In Section 5 we give our comparisons and in Section 6 we use our method to analyze three data sets.

## 2. Specifying the Bayes model

Given a PS distribution it is easy to write down the zero-inflated model (1). Specification of a Bayes model also requires a prior distribution. In this section, we present some forms and properties of (1) along with the Fisher information matrix that will be required for finding objective priors. We start with the PS distribution case and then specialize.

### 2.1. Zero-inflated power series distributions

Let $\mathbf{Y} = (Y_1, Y_2, \ldots, Y_n)$ be a random sample of size $n$, from $f^*(y|p,\theta)$ defined in (1), where $f(y|\theta)$ is given by

$$(5) \qquad\qquad f(y|\theta) = \frac{a_y \theta^y}{g(\theta)}, \qquad y = 0, 1, 2, \ldots,$$

in which $g(\theta) = \sum_{y=0}^{\infty} a_y \theta^y$ is the normalizing constant. It is easy to verify that $E(Y_1|p,\theta) = (1-p)\theta g'(\theta)/g(\theta)$. Writing

$$n_0 = \sum_{i=1}^{n} I_{[Y_i=0]}, \qquad S = \sum_{i=1}^{n} Y_i \quad \text{and} \quad \bar{Y} = S/n,$$

the likelihood function based on $\mathbf{Y}$ is

$$(6) \qquad\qquad L(p,\theta) = \{p + (1-p)f(0|\theta)\}^{n_0} \left(\frac{1-p}{g(\theta)}\right)^{n-n_0} \theta^S.$$

Using (6), it is an exercise to derive ML estimates for $(p,\theta)$. From (6) it is easy to derive that the per unit Fisher information matrix $I(p,\theta) = ((I_{ij}(p,\theta)))$ for $i, j = 1, 2$ is given by

$$I(p,\theta) =$$

$$\begin{pmatrix} \dfrac{1 - f(0|\theta)}{(1-p)\{p + (1-p)f(0|\theta)\}} & -\dfrac{\frac{g'(\theta)}{g(\theta)}f(0|\theta)}{p + (1-p)f(0|\theta)} \\[2em] -\dfrac{\frac{g'(\theta)}{g(\theta)}f(0|\theta)}{p + (1-p)f(0|\theta)} & (1-p)\left[-\dfrac{\left\{\frac{g'(\theta)}{g(\theta)}\right\}^2 (p + f(0|\theta))}{p + (1-p)f(0|\theta)} + \dfrac{1}{g(\theta)}\left(g''(\theta) + \dfrac{g'(\theta)}{\theta}\right)\right] \end{pmatrix}.$$

It is seen that the off-diagonal terms are nonzero. (In general, however, the off-diagonal terms are zero under the reparametrization $p^* = p + (1-p)f(0|\theta)$.)

Two special cases of (5) recur regularly, the zero-inflated Poisson and geometric. The zero-inflated binomial is similar; we do not treat it explicitly here.

### 2.1.1. Zero-inflated Poisson

The ZIP distribution with parameters $(p,\theta)$ results from (1) by using the Poisson $(\theta)$ probability mass function in place of $f(y|\theta)$ as indicated in (2). Parallel to (6), the likelihood function based on a sample of size $n$ is given by

$$(9) \qquad L(p, \theta) = \left\{ p + (1-p)e^{-\theta} \right\}^{n_0} \left\{ (1-p)e^{-\theta} \right\}^{(n-n_0)} \theta^s.$$

Using (9), the MLE for $(p, \theta)$ can be derived, see El-Shaarawi [10], as

$$(10) \qquad \hat{\theta}_1 = \frac{S(1 - e^{-\hat{\theta}_1})}{n - n_0} \quad \text{and} \quad \hat{p}_1 = \frac{\frac{n_0}{n} - e^{-\hat{\theta}_1}}{1 - e^{-\hat{\theta}_1}}.$$

Likewise, the test statistic for the score test for $\mathcal{H}_0 : p = 0$ can be derived as

$$(11a) \qquad T_s(\mathbf{Y}) = \frac{\left( \dfrac{n_0}{e^{-\hat{\theta}_0}} - n \right)^2}{n \left[ \dfrac{1 - e^{-\hat{\theta}_0}}{e^{-\hat{\theta}_0}} - \hat{\theta}_0 \right]},$$

where $\hat{\theta}_0 = \bar{Y}$ is MLE under $\mathcal{H}_0$, see Broek [1]. It can be shown that $\text{sgn}(\hat{p})\sqrt{T_s(\mathbf{Y})}$ asymptotically follows a standard normal distribution under $\mathcal{H}_0 : p = 0$ and a level $\alpha$ rejection region for testing $\mathcal{H}_0 : p = 0$ vs. $\mathcal{H}_1 : p > 0$ is given as

$$(11b) \qquad \text{sgn}(\hat{p})\sqrt{T_s(\mathbf{Y})} > z_\alpha,$$

where $z_\alpha$ is the upper $\alpha$ cut-off point from the standard normal distribution.

Similarly, the likelihood ratio test can be derived. We omit the details since, unlike the score test statistic, the likelihood ratio statistic does not have an explicit expression. If we denote the likelihood ratio test statistic by $T_l(\mathbf{Y})$, then $[\text{sgn}(\hat{p})\sqrt{T_l(\mathbf{Y})}]$ asymptotically follows $N(0, 1)$ under $\mathcal{H}_0 : p = 0$ and a level $\alpha$ rejection region for testing $\mathcal{H}_0 : p = 0$ vs. $\mathcal{H}_1 : p > 0$ is given as

$$(12) \qquad \text{sgn}(\hat{p})\sqrt{T_l(\mathbf{Y})} > z_\alpha.$$

For the ZIP family, the Fisher information matrix has entries

$$I(p, \theta) = \begin{pmatrix} \dfrac{1 - e^{-\theta}}{(1-p)\left\{ p + (1-p)e^{-\theta} \right\}} & -\dfrac{e^{-\theta}}{p + (1-p)e^{-\theta}} \\[2ex] -\dfrac{e^{-\theta}}{p + (1-p)e^{-\theta}} & \dfrac{1-p}{\theta} - \dfrac{p(1-p)e^{-\theta}}{p + (1-p)e^{-\theta}} \end{pmatrix}.$$

### 2.1.2. *Zero-inflated geometric*

The zero-inflated geometric distribution with parameters $(p, \theta)$ results from (1) by using the geometric $(\theta)$ probability mass function in place of $f(y|\theta)$ as indicated in (3). Parallel to (6), the likelihood function based on a sample of size $n$ is given by

$$(13) \qquad L(p, \theta) = \left\{ p + (1-p)(1-\theta) \right\}^{n_0} \left\{ (1-\theta)(1-p) \right\}^{n-n_0} \theta^s.$$

Using (13) the MLE's for $(p, \theta)$ can be derived; the test statistic for the score test for $\mathcal{H}_0$ is

$$T_s(Y) = \frac{n(1 + \bar{Y})}{\bar{Y}^2} \left[ \frac{n_0}{n}(1 + \bar{Y}) - 1 \right]^2.$$

In this case, the Fisher information matrix has entries

$$I(p,\theta) =$$
$$\begin{pmatrix} \dfrac{\theta}{(1-p)\left\{p+(1-p)(1-\theta)\right\}} & \dfrac{-1}{p+(1-p)(1-\theta)} \\[4mm] \dfrac{-1}{p+(1-p)(1-\theta)} & (1-p)\left\{\dfrac{\theta+(1-\theta)^2}{(1-\theta)^2\theta}+\dfrac{1-p}{p+(1-p)(1-\theta)}\right\} \end{pmatrix}.$$

For testing $\mathcal{H}_0 : p = 0$ vs. $\mathcal{H}_1 : p > 0$, the score test and likehood ratio test can be expressed in terms of rejection regions similar to those given in (11b) and (12).

## 2.2. Prior specification

It is well known that Jeffreys' prior is the reference prior in the absence of nuisance parameters, see Clarke and Barron [3]. That is, Jeffreys' prior is objective in the sense that using Jeffreys' prior gives a posterior that updates the prior as much as possible on average in relative entropy. Informally, it permits maximal information gain in a data transmission sense. For a small number of parameters, here 2, this is a reasonable optimality criterion.

By definition, Jeffreys' prior on $(p,\theta)$ is the square root of the determinant of the Fisher information matrix,

$$(14) \qquad\qquad \pi_J(p,\theta) \propto \left(\det(I(p,\theta))\right)^{1/2}.$$

For a zero-inflated power series distribution, there is no convenient expression in general for $\det(I(p,\theta))$. However, for the ZIP model, Jeffreys' prior is

$$(15\text{a}) \qquad\qquad \pi_J(p,\theta) \propto \frac{(1-e^{-\theta}-\theta e^{-\theta})^{1/2}}{[\theta\{p+(1-p)e^{-\theta}\}]^{1/2}}.$$

As is typical for reference priors, (15a) is improper: The integral over $\theta \in [0,\infty)$ diverges. Likewise, in a zero-inflated geometric, Jeffreys' prior is

$$(15\text{b}) \qquad\qquad \pi_J(p,\theta) \propto \frac{\theta^{1/2}}{(1-\theta)\{p+(1-p)(1-\theta)\}^{1/2}}.$$

Again, this is improper because the integral over $\theta$ diverges.

Jeffreys' prior, given by (14), would be appropriate if both $p$ and $\theta$ were of equal interest. Here, we are mainly interested in $p$. So, we used the Jeffreys' prior for $p$ for given $\theta$, that is

$$(16) \qquad\qquad \pi_J^c(p|\theta) \propto [I_{11}(p,\theta)]^{1/2},$$

and used the Jeffreys' prior for $\theta$ derived from the non-inflated model $f(y|\theta)$.

For a zero-inflated PS model, (16) gives

$$(17) \qquad\qquad \pi_J^c(p|\theta) = \frac{(1-f(0|\theta))^{1/2}}{\pi(1-p)^{1/2}\left[p+(1-p)f(0|\theta)\right]^{1/2}}.$$

If $g(\theta)$ corresponds to a zero-inflated Poisson model, (17) gives

$$(18\text{a}) \qquad\qquad \pi_J^c(p|\theta) = \frac{1}{\pi}\cdot\left[\frac{1-e^{-\theta}}{(1-p)\{p+(1-p)e^{-\theta}\}}\right]^{1/2},$$

and if $g(\theta)$ corresponds to a zero-inflated geometric distribution, (17) gives

$$(18b) \qquad \pi_J^c(p|\theta) = \frac{1}{\pi} \cdot \left[ \frac{\theta}{(1-p)\{p+(1-p)(1-\theta)\}} \right]^{1/2}.$$

Note that in both (18a,b) the range of $p$ includes a range of negative values, depending on $\theta$, for which the prior density is well defined.

Parallel to (16), the Jeffreys' prior for $\theta$ in the PS model is

$$(19) \qquad \pi_J(\theta) \propto \left[ \frac{g''(\theta)}{g(\theta)} + \frac{g'(\theta)}{\theta g(\theta)} - \left\{ \frac{g'(\theta)}{g(\theta)} \right\}^2 \right]^{1/2}.$$

Expression (19) gives $1/\sqrt{\theta}$ and $1/[(1-\theta)\sqrt{\theta}]$ for the Poisson and geometric model, respectively. These are improper. However, the posterior turns out to be proper because a finite number of data points suffice to make it so. If a proper joint objective prior is desired, Rissanen's prior, see Rissanen [16], can be adapted and gives similar results but is computationally more demanding.

## 3. Test criterion based on posterior probability

For the general case of a zero-inflated PS model, the posterior is formed by using (6) and (14). In the ZIP model, these expressions become (9) and (15a). Another reasonable choice would be (9) with (17), which becomes (18a), and $\pi(\theta) = 1/\sqrt{\theta}$.

Given these choices, the Bayes test for zero inflation is based on the posterior probability that $p > 0$. Thus, consider the statistic

$$(20) \qquad T(\mathbf{Y}) = P(p > 0|\mathbf{Y}) = \frac{\int_\Theta \int_0^1 L(u,\theta)\pi(u|\theta)\pi(\theta)dud\theta}{\int_\Theta \int_{\frac{-f(0|\theta)}{1-f(0|\theta)}}^1 L(u,\theta)\pi(u|\theta)\pi(\theta)dud\theta}.$$

The main point of this section is to derive an asymptotic test based on $T$ by finding the asymptotic distribution of $T(\mathbf{Y})$ under $\mathcal{H}_0 : p = 0$. It is reasonable to conclude that there is zero inflation when $P(p > 0|\mathbf{Y})$ is close to one. Consequently, from a frequentist standpoint, the rejection region is given by $T(\mathbf{Y}) > c$ where $c$ is chosen based on the given level of significance.

### 3.1. Finite sample properties of the test statistic

Note that (20) exploits the extended parameter space for $p$, namely, $-f(0|\theta)/(1-f(0|\theta)) < p < 1$ so that as $\theta$ increases, the lower bound approaches 0 from the left.

Let $P_{(p,\theta)}(\cdot)$ be the probability measure for a zero-inflated PS family. It can be verified that for large sample size, $P_{(p,\theta)}(T(\mathbf{Y}) > c)$ is increasing in $p$ for fixed $\theta$ and increasing in $\theta$ for fixed $p$. That is, as zero inflation increases the probability that $T$ is large (close to one) increases and that as the probability of large outcomes of $Y$ increases the probability of zero inflation also increases. This means that a single occurrence of zero can appear to be zero inflation if $\theta$ is large enough.

One feature which makes $T$ easy to use is that as a generality the joint posterior distribution for $(p, \theta)$ and the marginal posterior for $p$ are typically unimodal. Indeed, $\pi(p|\mathbf{y})$ is typically unimodal, even for small sample sizes. The posterior densities from the simulations reported in Section 5 and the data analysis in Section 6 are all unimodal.

### 3.2. *Asymptotic distribution of the test statistic*

For ease of exposition, let $\eta = (\eta_1, \eta_2) = (p, \theta)$ so that $Y_1, \ldots, Y_n$ are IID with density $f^*(\cdot|\eta)$ and $T(\mathbf{Y}) = P_\eta(\eta_1 > \eta_{10}|\mathbf{Y})$ where $\eta_{10}$ is a fixed value.

First, we sketch a proof that under $\eta = (\eta_{10}, \eta_{20}) = \eta_0$, the frequentist distribution of $T$ is asymptotically Uniform[0, 1], i.e., $T = U + o_p(1)$ as $n \to \infty$. Then we derive an expression for the asymptotic behavior of the first two moments of $T$. Although these arguments are presented in the zero-inflated PS context, they appear to be more general.

Start by writing

$$\ell(\eta) = (1/n) \sum_{i=1}^n \log f^*(y_i|\eta) \quad \text{and} \quad \hat{\eta} = (\hat{\eta}_1, \hat{\eta}_2) = \arg\max_\eta \ell(\eta).$$

Letting $D_j$ denote $\partial/\partial\eta_j$ for $j = 1, 2$ define

$$a_{ij} = D_i D_j \ell(\hat{\eta}) \quad \text{and} \quad a_{ijk} = D_i D_j D_k \ell(\hat{\eta}).$$

Under consistency conditions for the MLE and expected local supremum conditions on $f^*(\cdot|\eta)$ on a neighborhood around a fixed value $\eta_0$,

$$a_{ijk} \to E_{\eta_0} D_i D_j D_k \log f_\eta^*(Y_1|\eta_0), \quad a.e., P_{\eta_0}.$$

The empirical Fisher information is $I(\hat{\eta}) = (I_{ij}(\hat{\eta}))$ and it is seen that $I_{ij}(\hat{\eta}) = -(1/n)a_{ij}(\hat{\eta})$. To ensure $I(\hat{\eta})$ is well defined, assume that it is positive definite on a set $S^*$ with $P_\eta(S^*) = 1 + o(1/\sqrt{n})$. Now the inverse is $I^{-1}(\hat{\eta}) = (I^{ij}(\hat{\eta}))$; it is needed to define the quantities that will appear in the asymptotic expression for $T$. Set

$$m_i(\hat{\eta}) = \frac{I^{i1}(\hat{\eta})}{I^{11}(\hat{\eta})}, K^{ij}(\hat{\eta}) = I^{ij}(\hat{\eta}) - \frac{I^{i1}(\hat{\eta})I^{j1}(\hat{\eta})}{I^{11}(\hat{\eta})}$$

and denote $\hat{\pi}_j(\hat{\eta}) = D_j\pi(\hat{\eta})$. Finally, the quantities that appear in the asymptotic expression to order $O(1/\sqrt{n})$ are the second degree Hermite polynomial $J_2(t) = t^2 - 1$, and two correction terms

$$G_3(\hat{\eta}) = \frac{1}{6}a_{ijk}m_i m_j m_k (I^{11}(\hat{\eta}))^{3/2}$$

and

$$G_1(\pi, \hat{\eta}) = \frac{\hat{\pi}_j m_j}{\pi(\hat{\eta})}\sqrt{I^{11}} + \frac{1}{2}a_{ijk}K^{ij}m_k\sqrt{I^{11}} + \frac{1}{2}a_{ijk}m_i m_j m_k (I^{11}(\hat{\eta}))^{3/2},$$

using the convention that repeated indices indicate summation. To get the form of the result, let

$$W = \sqrt{\frac{n}{I^{11}(\hat{\eta})}}(\eta_{10} - \hat{\eta}_1) \quad \text{and} \quad V = \sqrt{\frac{n}{I^{11}(\hat{\eta})}}(\eta_1 - \hat{\eta}_1).$$

Note that under $\eta_1 = \eta_{10}$, $V$ is the same as $W$.

At last, from (2.3.19) in Datta and Mukerjee [6], taking $\beta_1 = \beta_2 = 0$ we get

$$\begin{aligned}
P(\eta_1 &\leq \eta_{10}|\mathbf{Y}) \\
&= P(v \leq w|\mathbf{Y}) \\
(21) \quad &= \Phi(w) + n^{-1/2}\phi(w)\{G_1(\pi, \hat{\eta}) + G_3(\hat{\eta})J_2(w)\} + o_p(n^{-1/2}),
\end{aligned}$$

where $w$ is the observed value of $W$ and $\phi(\cdot)$ and $\Phi(\cdot)$ are standard normal pdf and cdf, respectively. However, when $\eta_{10}$ is true $W$ is asymptotically $N(0, 1)$. So, by the inverse probability integral transform $\Phi(W)$ is Uniform[0, 1] and the $1/\sqrt{n}$ terms ensure the required rate. Thus, since $T$ is of the form $1 - P(\eta_1 \leq \eta_{10}|\mathbf{Y})$, $T$ is asymptotically Uniform[0, 1] as well.

To derive expressions for the moments of $T$, write $G_1(\pi, \hat\eta) = \Gamma_1(\eta_0) + o_p(1)$, and $G_3(\hat\eta) = \Gamma(\eta_0) + o(1)$. Recognizing that $J_2$ is just a polynomial, it can be seen that there is an $H(\eta_0)$ so that (21) can be written as

$$(22) \qquad P(\eta_1 \leq \eta_{10}|\mathbf{Y}) = P(V \leq w|\mathbf{Y}) = \Phi(w) + n^{-1/2}\phi(w)H(\eta_0) + o_p(n^{-1/2}),$$

and the expectation with respect to $P_\eta$ can be taken on both sides. Using the result from that and applying Step 3 from Datta and Mukerjee [6], page 19, gives an expression for the frequentist probability from the middle term in (22):

$$(23) \qquad P_{\eta_0}(V \leq w) = \Phi(w) + n^{-1/2}\phi(w)H^*(\eta_0) + o_p(n^{-1/2}),$$

where $H^*$ is derived from $H$ and the $\eta_1$ in the probability is $\eta_{10}$, the true value. Differentiating (23) it is possible to derive an approximation for the density of $V$, $f_V(w|\eta_0)$.

Finally, by using (21) and $f_W(w|\eta_0)$, it is possible to derive expressions for the first 2 moments of $T$, because they only depend on $W$. Doing so gives that they are $1/2$ and $1/12$, as expected from the limiting uniform. However, given the $1/\sqrt{n}$ correction terms, it is possible to equate the expressions for the first 2 moments of $P(\eta_1 \leq \eta_{10}|\mathbf{Y})$ to the first two moments of a $Beta(\alpha, \beta)$ and thereby derive expressions for $\alpha$ and $\beta$. Obviously, the resulting $\hat\alpha$ and $\hat\beta$ must converge to 1, i.e., give the Uniform[0, 1] in the limit for large $n$, but for finite $n$ this provides a more refined approximation.

## 4. Credible intervals

Although one can in principle find a Bayes estimate for $p$, under say squared error loss, and find its posterior variance, Bayes tests are based on posterior probabilities which in turn are based on the posterior density. These also lead to credible sets.

There are two main types of credible sets. The first is analogous to confidence intervals: $\alpha/2$ of the probability in the tails is clipped off and the upper and lower boundaries announced. The second is highest posterior density HPD, i.e., a set of the form $R(\pi_\alpha) = \{p : \pi(p|\mathbf{y}) \geq \pi_\alpha\}$, where $\pi_\alpha$ is the largest constant such that $P(p \in R(\pi_\alpha)|\mathbf{y}) \geq 1 - \alpha$. For symmetric unimodal densities the two types of interval are equal, and here HPD sets are obtained from a variant on the procedure used to get $\alpha$-clipped credible sets. The basic idea is that if the credible interval or HPD set for $p$ contains the value $p = 0$, then we may conclude that there is not enough evidence of zero inflation in the data.

Difficulties in the cases studied here arise because the posterior is not available in a convenient analytic form: The priors discussed in Section 2.2 do not yield tractable marginal posteriors for $p$ by directly integrating $\theta$ out of joint posterior. Consequently, we find a Markov chain Monte Carlo (MCMC) estimate of the marginal posterior distribution and use it to find the $1 - \alpha$ credible and HPD sets. Thus, given a sample from the marginal posterior $\pi(p|\mathbf{Y} = \mathbf{y})$ it is easy to form a $1 - \alpha$ credible interval by choosing the $\alpha/2$ and $1 - \alpha/2$ sample quantiles. This can also be done using draws from the joint $(p, \theta)$ posterior density. The HPD set can be found by using the draws to estimate $\pi(p|\mathbf{Y} = \mathbf{y})$ by, say, $\hat{\pi}(p|\mathbf{Y} = \mathbf{y})$ and obtaining approximate HPD sets from $\hat{\pi}$.

Suppose that $\pi(p, \theta|\mathbf{y})$ and $\pi(p|\mathbf{y})$ are the joint and marginal posterior, respectively, so that

$$\pi(p|\mathbf{y}) = \int_0^\infty \pi(p, \theta|\mathbf{y}) d\theta.$$

In the case of ZIP model (9), with the conditional Jeffreys' prior (18a) for $p$, and Jeffreys' prior $\pi(\theta) = 1/\sqrt{\theta}$ for $\theta$, the joint posterior is based on

$$(24) \qquad \pi(p, \theta|\mathbf{y}) \propto \left\{ p + (1-p)e^{-\theta} \right\}^{n_0 - 1/2} (1-p)^{n - n_0 - 1/2} e^{-\theta(n - n_0)} \theta^{s - 1/2}.$$

Using (24), the goal is to estimate $\pi(p|\mathbf{y})$ from a joint sample of $(p, \theta)$ drawn from $\pi(p, \theta|\mathbf{y})$. Let $\{(p^{(i)}, \theta^{(i)}), i = 1, \ldots, B\}$ be an MCMC sample from $\pi(p, \theta|\mathbf{y})$ so that $\pi(p|\mathbf{y})$ can be estimated at $p = p^{(j)}$ by

$$(25) \qquad \hat{\pi}(p^{(j)}|\mathbf{y}) = \frac{1}{B} \sum_{i=1}^{B} \pi(p^{(j)}, \theta^{(i)}|\mathbf{y}).$$

Since it is computationally difficult to draw samples $(p^{(i)}, \theta^{(i)})$ from (24), we use a reparametrized model by transforming $p^* = p + (1-p)e^{-\theta}$. Incidentally, note that the parameters $p^*$ and $\theta$ result in an orthogonal reparameterization of the ZIP model. It can be checked that the Fisher information matrix is diagonal given by

$$I(p^*, \theta) = \text{diag}\left( \frac{1}{p^*(1-p^*)}, \frac{(1 - e^{-\theta} - \theta e^{-\theta})(1-p^*)}{\theta(1 - e^{-\theta})^2} \right).$$

As a result of this reparameterization, the joint posterior for $(p^*, \theta)$ can written as a product of their marginals. In fact this idea can be extended in general for all zero-inflated PS distributions.

Therefore, using the above fact the joint posterior distribution can be written as

$$(26) \qquad \pi(p^*, \theta|\mathbf{y}) \propto (p^*)^{n_0 - 1/2} (1-p^*)^{n - n_0 - 1/2} \left( \frac{e^{-\theta}}{1 - e^{-\theta}} \right)^{n - n_0} \theta^{s - 1/2}.$$

From (26) it is seen that $(p^*|\mathbf{y})$ follows a Beta$(n_0 + 1/2, n - n_0 + 1/2)$, so it is easy to draw posterior samples of $p^*$. To draw samples of $\theta$, we use rejection sampling with a suitably chosen gamma distribution as envelope. In this way it is possible to generate a representative sample $\{(p^{*(i)}, \theta^{(i)}), i = 1, \ldots, B\}$ from the joint posterior and using the relationship between $p^*$ and $p$, we get $\{(p^{(i)}, \theta^{(i)}), i = 1, \ldots, B\}$ where $p^{(i)} = (p^{*(i)} - e^{-\theta^{(i)}})/(1 - e^{-\theta^{(i)}})$. Subsequently, using (25), we get an estimate of the marginal posterior density $\pi(p|\mathbf{y})$ at $p = p^{(j)}$.

A $100(1-\alpha)\%$ Bayesian credible interval for $p$ is simply $(p_{(\alpha/2)}, p_{(1-\alpha/2)})$, where $p_{(k)}$ is the $k$-th quantile of $\{p^{(i)}, i = 1, \ldots, B\}$. To find the HPD interval, there are several methods and algorithms, see Chen and Shao [2], and the references therein. Here, using the HPD set from the posterior sample $\{p^{(i)}, i = 1, \ldots, B\}$, we find $\{\hat{\pi}(p^{(i)}|\mathbf{y}), i = 1, \ldots, B\}$ and set $\pi_\alpha$ to be the $100\alpha$-th percentile of $\hat{\pi}(p^{(i)}|\mathbf{y})$. This is adequate because the estimated posteriors are unimodal. Once we get $\pi_\alpha$, we solve $\hat{\pi}(p|\mathbf{y}) = \pi_\alpha$ for cut-off values of $p$ to find the lower and upper limits of the HPD interval.

## 5. Performance comparison

In this section the performance of the test statistic (20) for $\mathcal{H}_0 : p = 0$ vs. $\mathcal{H}_1 : p > 0$ in the ZIP family is compared to the performance of the score test and the likelihood ratio test in simulations. Recall that the Bayes test is formed from (9), (18a), and $\pi(\theta) = 1/\sqrt{\theta}$. The score test is given by (11a) and we numerically obtain the likelihood ratio statistic. In the Table 1 below, we have computed the power of these three tests for several choices of $n$, $p$, $\theta$ for $\mathcal{H}_0$ vs $\mathcal{H}_1$ for level $\alpha = 0.05$; in Table 2 the power of the one-sided Bayes test is compared to the power of the two-sided score and LR tests as well, also at the $\alpha = 0.05$ level. We can see that the Bayes test performs somewhat better than the two-sided score test and the two-sided likelihood ratio test.

However, as in Section 5, the simulations for the Bayes test require MCMC because calculating $P(p > 0|\mathbf{Y})$ under the ZIP model is not straightforward. Indeed, in general, it is not possible to provide a procedure that will work for any zero-inflated PS distribution. Nevertheless, for the ZIP model, (26), implies that (20)

TABLE 1

*The entries are the powers for the Bayesian, one-sided score, and one-sided LR tests, with 10,000 simulations. The asterisks indicate when the values where the Bayes test has highest power. They are clustered around small to moderate p, small θ and small to moderate n*

| | $p$ | 0.00 | | | 0.10 | | | 0.30 | | | 0.40 | | |
|---|---|---|---|---|---|---|---|---|---|---|---|---|---|
| $\theta$ | $n$ | Score | Bayes | LR | Score | Bayes | LR | Score | Bayes | LR | Score | Bayes | LR |
| 0.5 | 20 | 0.049 | 0.045 | 0.047 | 0.065 | 0.068 | 0.064 | 0.111 | 0.105 | 0.103 | 0.144 | 0.134 | 0.118 |
| | 50 | 0.046 | 0.043 | 0.042 | 0.078 | 0.076 | 0.072 | 0.180 | 0.159 | 0.154 | 0.251 | 0.212 | 0.209 |
| | 100 | 0.050 | 0.047 | 0.046 | 0.096 | 0.090 | 0.081 | 0.284 | 0.262 | 0.263 | 0.376 | 0.363 | 0.345 |
| 1.0 | 20 | 0.040 | 0.049 | 0.036 | 0.083 | 0.094* | 0.082 | 0.232 | 0.247* | 0.228 | 0.318 | 0.323* | 0.311 |
| | 50 | 0.040 | 0.049 | 0.040 | 0.123 | 0.133* | 0.126 | 0.433 | 0.434* | 0.417 | 0.585 | 0.582 | 0.566 |
| | 100 | 0.045 | 0.047 | 0.048 | 0.181 | 0.182 | 0.188 | 0.670 | 0.671 | 0.680 | 0.840 | 0.841 | 0.843 |
| 1.5 | 20 | 0.042 | 0.053 | 0.040 | 0.123 | 0.143* | 0.116 | 0.389 | 0.420* | 0.387 | 0.544 | 0.564* | 0.537 |
| | 50 | 0.040 | 0.047 | 0.043 | 0.214 | 0.225* | 0.212 | 0.730 | 0.747* | 0.739 | 0.884 | 0.895* | 0.888 |
| | 100 | 0.045 | 0.046 | 0.046 | 0.345 | 0.311 | 0.351 | 0.951 | 0.936 | 0.953 | 0.992 | 0.991 | 0.993 |
| 2.0 | 20 | 0.046 | 0.052 | 0.035 | 0.194 | 0.213* | 0.175 | 0.615 | 0.649* | 0.600 | 0.763 | 0.801* | 0.758 |
| | 50 | 0.053 | 0.053 | 0.045 | 0.345 | 0.363* | 0.346 | 0.936 | 0.93 | 0.935 | 0.988 | 0.988 | 0.986 |
| | 100 | 0.044 | 0.053 | 0.042 | 0.577 | 0.484 | 0.557 | 0.998 | 0.995 | 0.998 | 1.000 | 1.000 | 1.000 |

TABLE 2

*The entries are the powers for the Bayesian test and the two-sided score and two-sided LR tests, with 10,000 simulations. Putting asterisks in this table gives the same pattern as in Table 1, but stronger*

| | $p$ | 0.00 | | | 0.10 | | | 0.30 | | | 0.40 | | |
|---|---|---|---|---|---|---|---|---|---|---|---|---|---|
| $\theta$ | $n$ | Score | Bayes | LR | Score | Bayes | LR | Score | Bayes | LR | Score | Bayes | LR |
| 0.5 | 20 | 0.045 | 0.045 | 0.061 | 0.043 | 0.068 | 0.052 | 0.065 | 0.105 | 0.057 | 0.087 | 0.134 | 0.066 |
| | 50 | 0.046 | 0.043 | 0.050 | 0.055 | 0.076 | 0.056 | 0.122 | 0.159 | 0.106 | 0.181 | 0.212 | 0.136 |
| | 100 | 0.051 | 0.047 | 0.051 | 0.066 | 0.090 | 0.058 | 0.185 | 0.262 | 0.174 | 0.277 | 0.363 | 0.248 |
| 1.0 | 20 | 0.048 | 0.049 | 0.058 | 0.057 | 0.094 | 0.062 | 0.142 | 0.247 | 0.143 | 0.203 | 0.323 | 0.198 |
| | 50 | 0.049 | 0.049 | 0.051 | 0.075 | 0.133 | 0.078 | 0.303 | 0.434 | 0.296 | 0.443 | 0.582 | 0.430 |
| | 100 | 0.052 | 0.047 | 0.050 | 0.117 | 0.182 | 0.115 | 0.571 | 0.671 | 0.542 | 0.767 | 0.841 | 0.739 |
| 1.5 | 20 | 0.047 | 0.053 | 0.057 | 0.081 | 0.143 | 0.074 | 0.280 | 0.420 | 0.267 | 0.411 | 0.564 | 0.409 |
| | 50 | 0.051 | 0.047 | 0.051 | 0.140 | 0.225 | 0.131 | 0.618 | 0.747 | 0.612 | 0.806 | 0.895 | 0.809 |
| | 100 | 0.049 | 0.046 | 0.054 | 0.244 | 0.311 | 0.236 | 0.913 | 0.936 | 0.908 | 0.983 | 0.991 | 0.982 |
| 2.0 | 20 | 0.041 | 0.052 | 0.071 | 0.128 | 0.213 | 0.113 | 0.501 | 0.649 | 0.471 | 0.670 | 0.801 | 0.644 |
| | 50 | 0.049 | 0.053 | 0.057 | 0.257 | 0.363 | 0.228 | 0.890 | 0.935 | 0.880 | 0.975 | 0.988 | 0.973 |
| | 100 | 0.047 | 0.053 | 0.045 | 0.451 | 0.484 | 0.440 | 0.995 | 0.995 | 0.994 | 1.000 | 1.000 | 1.000 |

can be written as

$$P(p > 0|\mathbf{Y})$$

$$= \frac{\int_0^\infty \int_0^1 w(p^*, \theta)(1-p^*)^{n-n_0} e^{-\theta(n-n_0)} \theta^{s-(n-n_0)} \pi^*(p^*, \theta) dp^* d\theta}{\int_0^\infty \int_0^1 \left(\frac{\theta}{1-e^{-\theta}}\right)^{n-n_0} p^{*n_0}(1-p^*)^{n-n_0} e^{-\theta(n-n_0)} \theta^{s-(n-n_0)} \pi^*(p^*, \theta) dp^* d\theta}$$

$$(27) \quad = \frac{E^g\left[w(p^*, \theta)\pi^*(p^*, \theta)\right]}{E^g\left[\left(\frac{\theta}{1-e^{-\theta}}\right)^{n-n_0} \pi^*(p^*, \theta)\right]}$$

where $\pi^*(p^*, \theta)$ is the joint prior of $(p^*, \theta)$,

$$g(p^*, \theta) = p^{*n_0}(1-p^*)^{n-n_0} e^{-\theta(n-n_0)} \theta^{s-(n-n_0)},$$

and $w(p^*, \theta) = I_{[p^* > e^{-\theta}]} \left(\frac{\theta}{1-e^{-\theta}}\right)^{n-n_0}$.

So, we draw a random sample $\{(p^{*(i)}, \theta^{(i)}), i = 1, \ldots, B\}$ where $p^{*(i)} \sim Beta(n_0 + 1, n - n_0 + 1)$ and $\theta^{(i)} \sim gamma(n - n_0, s - (n - n_0) + 1))$ and calculate

$$P(p > 0|\mathbf{Y}) = \frac{\sum_{i=1}^B I_{[p^{*(i)} > e^{-\theta^{(i)}}]} \pi^*(p^{*(i)}, \theta^{(i)}) \left(\frac{\theta^{(i)}}{1-e^{-\theta^{(i)}}}\right)^{n-n_0}}{\sum_{i=1}^B \pi^*(p^{*(i)}, \theta^{(i)}) \left(\frac{\theta^{(i)}}{1-e^{-\theta^{(i)}}}\right)^{n-n_0}}.$$

In Tables 1 and 2, the test statistic $T(\mathbf{Y})$ from (20) is compared to the cut-off point found from the asymptotic distribution of $T(\mathbf{Y})$ under $\mathcal{H}_0$, i.e., we use upper $\alpha$ point on $uniform(0, 1)$ as described in Section 3.2. All the simulations are based on 10,000 replications with $B = 10,000$ MCMC samples in each replication.

We comment that for the case of a zero-inflated geometric distribution, the procedure is a little easier: It is just a matter of drawing samples from two different Beta distributions with parameters based on the sample. So, it is easy to find an MCMC estimate of the test statistic.

Table 1 shows that for the one-sided test, all 3 tests have roughly the same level when $p = 0$. In fairness, the level for the score and LR tests is a little lower leading to lower power against alternatives. However, looking at how the power of all three tests indicated rises as $p$ increases, it is clear that for mid-sized $p$ and smallish $\theta$ the Bayes test has noticeably higher power, especially for small $n$. In fact, a good test is most important on this range because it is hard to distinguish zero inflation from its absence when $\theta$ is small or moderate, $p$ ranges from small to mid-sized values, and $n$ is not large.

Table 2 is similar to Table 1 but the score and LR tests are two-sided. It shows that the same properties hold, but a little more strongly. This may be attributed to the fact that the Bayes test uses an extended parameter space, allowing some mass to represent zero deflation.

## 6. Data analysis

To demonstrate the efficacy of our technique, we apply it to test for presence of zero inflation in three famous datasets. We also give comparative values from other techniques. In general, the results from the techniques corroborate each other so the fact they are based on different principles lends credence to the conclusions.

The first dataset that we look at is the Urinary Tract Infection (UTI) data used in Broek [1] who used a score test to detect zero inflation in a Poisson model. The data are collected from 98 HIV-infected men, attending the department of internal medicine at the Utrecht University hospital. The number of times they had a urinary tract infection was recorded as $Y$. The data are recorded in Table 3. Merely by looking at the data it is clear that zero inflation is present.

Our method yields a Bayes factor for testing $\mathcal{H}_0 : p = 0$ vs. $\mathcal{H}_1 : p > 0$ of $B_{10} = 223.13$. The details of computation of Bayes factor will be reported elsewhere. This is strong evidence in favor of the alternative, which is no surprise. In fact, $P(p > 0|\mathbf{y}) = .999$. The observed value of the score statistic is 15.34 giving a $p$-value 0.0001.

The next data set we consider is the Terrorism data from Conigliani, Castro and O'Hagan [5]. Table 4 shows the data concerning the number of incidents of international terrorism per month ($Y$) in the United States between 1968 and 1974. It is not immediately clear if there is a zero-inflation in this data set. Conigliani, Castro and O'Hagan [5] find a Fractional Bayes factor for this data set of 0.0089; we find a Bayes factor of $B_{10} = 0.28$. In fact, $P(p > 0|\mathbf{y}) = 0.507$, an indeterminate value. The observed value of the score statistic is 0.04 and a $p$-value 0.83. All three assessments agree that there is no evidence of zero inflation.

The third data set we analyzed is the Cholera data first analyzed by McKendrick [15]. Table 5 shows the number of patients per household suffering from cholera in a village in India in 1920's. Again, looking at the data strongly suggests zero inflation. While the Bayes factor is $B_{10} = 238090$, very strong evidence for zero inflation, under our method, $P(p > 0|\mathbf{y}) = .9999$. The observed value of the score statistic is 30.56, effectively giving a $p$-value of 0. Again, all three assessments agree for this example.

Although tests are useful for quantifying degree of belief, they are not the same as looking at the posterior distributions directly. Figure 1 shows plots of the marginal posteriors for $p$ resulting from applying the ZIP model to each of the three data sets. All three posteriors are unimodal and appear roughly symmetric. The location of the mode, and the spread around it determine the most credible values of $p$. For the UTI and Cholera data the determination is clear: Substantial zero inflation is present. For the Terror data, the graph does not give a clear answer. The slight asymmetry makes it difficult to tell whether $p = 0$ is reasonable. In fact, the test shows it is, but this would be open to question from merely looking at the diagram.

Table 6 gives 95% Bayesian credible and HPD intervals for the three data sets under consideration. Also the marginal posterior distributions of $p$ are given in

TABLE 3
*UTI data*

| $Y$ | 0 | 1 | 2 | 3 | Total |
|---|---|---|---|---|---|
| Frequency | 81 | 9 | 7 | 1 | 98 |

TABLE 4
*Terror data*

| $Y$ | 0 | 1 | 2 | 3 | 4 | Total |
|---|---|---|---|---|---|---|
| Frequency | 38 | 26 | 8 | 2 | 1 | 75 |

TABLE 5
*Cholera data*

| $Y$ | 0 | 1 | 2 | 3 | 4 | Total |
|---|---|---|---|---|---|---|
| Frequency | 168 | 32 | 16 | 6 | 1 | 223 |

TABLE 6
*Bayesian credible and HPD intervals*

| Data | Credible Interval | HPD Interval |
|------|-------------------|--------------|
| Terror | $(-0.6735, 0.2945)$ | $(-0.5560, 0.3654)$ |
| Cholera | $(0.4619, 0.7095)$ | $(0.4700, 0.7144)$ |
| UTI | $(0.3433, 0.8240)$ | $(0.4271, 0.8561)$ |

Figure 1. From the intervals and the figures as well, it is evident that there is noticeable amount of zero inflation in Cholera data and UTI data because the interval of concentration of the posterior distributions does not include zero whereas for Terror data the posterior distribution of $p$ is centered around zero and the interval contains zero, signifying the absence of zero inflation in the data.

Finally, for the sake of completeness, Table 6 gives 0.95 credible intervals and HPD intervals calculated from the posteriors. It is seen that for the Cholera and UTI data that 0 is not in the intervals. This is consistent with the presence of zero inflation. For the Terror data, 0 is in the interval. The interval is so wide much of it includes negative values. So, it is not a surprise that zero inflation is not indicated by the test. Note that the credible and HPD sets are close for the cholera data indicating symmetry, but for the other two data sets the difference in the intervals suggests some left skewing, more for Terror than for Cholera.

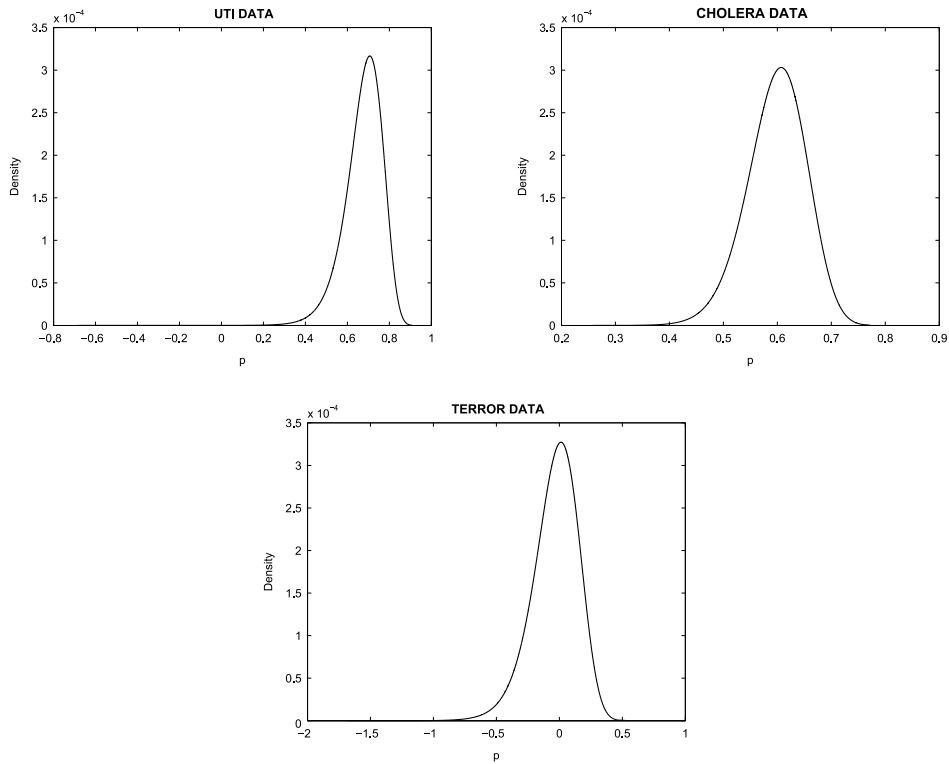

Fig 1. *Estimated posterior densities of p.*

## 7. Conclusions

Overall, this article gives a general Bayesian setup for testing for zero inflation in PS distributions that can be compared to existing likelihood based methods occurring in frequentist treatments. The basic idea is to extend the parameter space to include a small range of negative values for the weight on zero inflation. Thus, the null hypothesis $\mathcal{H}_0 : p = 0$ becomes an interior point of the parameter space and a standard Bayesian approach is feasible.

Our simulations suggest the Bayesian test has power as high as, or slightly higher than the likelihood based tests, even when objective priors that are somewhat unfavorable to the hypothesis $\mathcal{H}_0 : p = 0$ are used to automate the procedure. Interval estimation for $p$ proceeds similarly, using the extended parameter space.

The technique of extending the parameter space applies generally to Bayes, and potentially to frequentist, testing for zero inflation with count data, but obviously can apply to many situations where two distributions are mixed and one wants to know whether one component can be set to zero. In fact, the asymptotics for this test require only generic regularity conditions; they do not rely on specific forms of the likelihood such as exponential families. A further test of the method, aside from applying it to more general mixtures, would be extending it to a class of regression problems by including covariates.

**Acknowledgments.** The authors would like to thank Prof. P. K. Sen and Prof. J. K. Ghosh for their valuable comments. They are also thankful to Professor Edsel Peña for his comments that lead to improvement in the presentation. The authors acknowledge Prof. Xiao-Li Meng for sharing the Cholera data. The authors thank Professor Tathagata Bandyopadhyay for drawing their attention to the ZIP problem. Part of this research was conducted while G. S. Datta was a SAMSI/Duke University visiting fellow and B. Clarke was a SAMSI visitor. The support of both the institutes is gratefully acknowledged. A. Bhattacharya received student support from this grant.

## References

[1] BROEK, J. V. D. (1995). A score test for zero inflation on a Poisson distribution. *Biometrics* **51** 738–743. MR1349912
[2] CHEN, M.-H. AND SHAO, Q.-M. (1999). Monte Carlo estimation of Bayesian credible and HPD intervals. *J. Comput. Graph. Statist.* **8** 69–92. MR1705909
[3] CLARKE, B. AND BARRON, A. (1994). Jeffreys prior is asymptotically least favorable under entropy loss. *J. Statist. Plann. Inference* **41** 37–60. MR1292146
[4] COCHRAN, W. G. (1954). Some methods of strengthening $\chi^2$ tests. *Biometrics* **10** 417–451. MR0067428
[5] CONIGLIANI, C., CASTRO, J. I. AND O'HAGAN, A. (2000). Bayesian assessment of goodness of fit against nonparametric alternatives. *Canad. J. Statist.* **28** 327–342. MR1792053
[6] DATTA, G. S. AND MUKERJEE, R. (2004). *Probability Matching Priors*: *Higher Order Asymptotics*. Springer, Berlin. MR2053794
[7] DENG, D. AND PAUL, S. R. (2000). Score test for zero inflation in generalized linear models. *Canadian J. Statist.* **28** 563–570. MR1793111
[8] DENG, D. AND PAUL, S. R. (2005). Score tests for zero-inflation and over-dispersion in generalized linear models. *Statist. Sinica* **15** 257–276. MR2125731

[9] Efron, B. (1986). Double exponential families and their use in generalized linear regression. *J. Amer. Statist. Assoc.* **81** 709–721. MR0860505

[10] El-Shaarawi, A. H. (1985). Some goodness-of-fit methods for the Poisson plus added zeros distribution. *Appl. Environ. Microbiology* **49** 1304–1306.

[11] Ghosh, S. K., Mukhopadhyay, P. and Lu, J. C. (2006). Bayesian analysis of zero-inflated regression models. *J. Statist. Plann. Inference* **136** 1360–1375. MR2253768

[12] Hall, D. B. (2000). Zero-inflated Poisson and binomial regression with random effects: A case study. *Biometrics* **56** 1030–1039. MR1815581

[13] Johnson, N. L., Kotz, S. and Kemp, A. W. (1992). *Univariate Discrete Distributions*, 2nd ed. Wiley, New York. MR1224449

[14] Lambert, D. (1992). Zero-inflated Poisson regression, with an application to defects in manufacturing. *Technometrics* **34** 1–14.

[15] McKendrick, A. G. (1926). Application of mathematics to medical problems. *Proc. Edin. Math. Soc* **44** 98–130.

[16] Rissanen, J. (1983). A universal prior for integers and estimation by minimum description length. *Ann. Statist.* **11** 416–431. MR0696056

[17] Rao, C. R. and Chakravarti, I. M. (1956). Some small sample tests of significance for a Poisson distribution. *Biometrics* **12** 264–282. MR0081596

[18] Self, S. G. and Liang, K. Y. (1987). Asymptotic properties of maximum likelihood estimators and likelihood ratio tests under nonstandard conditions. *J. Amer. Statist. Assoc.* **82** 605–610. MR0898365

[19] Silvapulle, M. S. and Silvapulle, P. (1995). A score test against one-sided alternatives. *J. Amer. Statist. Assoc.* **90** 342–349. MR1325141